\documentclass[11pt,a4paper]{article}
\usepackage[a4paper,margin=2.5cm]{geometry}
\usepackage[T1]{fontenc}
\usepackage{amsmath,amsthm}
\usepackage{mathptmx,amssymb,courier}
\usepackage{microtype}
\usepackage{enumitem}
\usepackage[round,authoryear]{natbib}
\usepackage[hidelinks]{hyperref}
\newtheorem{theorem}{Theorem}[section]
\newtheorem{proposition}[theorem]{Proposition}
\newtheorem{lemma}[theorem]{Lemma}
\newtheorem{corollary}[theorem]{Corollary}
\newtheorem{conjecture}[theorem]{Conjecture}
\theoremstyle{definition}
\newtheorem{definition}[theorem]{Definition}
\newtheorem{example}[theorem]{Example}
\theoremstyle{remark}
\newtheorem{remark}[theorem]{Remark}

\newcommand{\M}{\mathfrak M}
\newcommand{\Sym}{\mathbb S}
\newcommand{\tr}{\operatorname{tr}}
\newcommand{\Tan}{\operatorname{Tan}}
\newcommand{\qstr}{q_{\mathrm{str}}}
\newcommand{\cS}{\mathcal S}
\newcommand{\cL}{\mathcal L}
\newcommand{\cN}{\mathcal N}
\newcommand{\cB}{\mathfrak B}
\newcommand{\cE}{\mathfrak E}
\newcommand{\dbar}{\bar d}
\newcommand{\Dl}{\underline{\Delta}}
\newcommand{\Du}{\overline{\Delta}}

\title{Structural latent dimension:\\ deterministic geometry and local Fisher structure}
\author{Matteo Farn\`e\\[2pt]
{\small Department of Statistical Sciences, University of Bologna}\\
{\small\texttt{matteo.farne@unibo.it}}}
\date{November 2026}

\begin{document}
\maketitle

\begin{abstract}
We study the \emph{structural dimension} of a second-order object: the minimum latent dimension that remains necessary relative to a declared class of competing dependence. Building on minimum-rank factor analysis and low-rank-plus-sparse decomposition, we use distances to lower-dimensional classes to describe inference near structural ambiguity. General distance and tangent-cone facts yield confidence statements and local margins. For rank-zero competitors with bounded residual degree, we derive exact maximum-weight $b$-matching formulas, support-slack tangent cones and sparse-spike angles; a shared residual support across a field can produce separation absent pointwise. We give positive-semidefinite brackets, a dimension-free Gaussian likelihood expansion, an exact quotient by the tangent of a low-rank background, and a detection floor from the margin. Equal-angle configurations with different branch multiplicities motivate static and spectral detection conjectures. Joint low-rank/residual profiling and the sharp composite information boundaries remain open.
\end{abstract}


\section{Introduction}

Factor and latent-variable models represent high-dimensional second-order objects by a low-rank component and a structured remainder. Factor-number criteria relate dimension to signal strength under their model assumptions \citep{BaiNg2002}. Minimum-rank factor analysis studies the least common rank compatible with diagonal noise \citep{tenBerge1991,Shapiro2002,Saunderson2012}. Low-rank-plus-sparse decomposition jointly recovers components with unknown rank and support, using geometric conditions for identifiability and exact recovery \citep{Chandrasekaran2011,Chandrasekaran2012}. These theories already recognise that a low-rank direction can also have a sparse explanation. Strength and necessity can then diverge: a direction with a large eigenvalue may be reproducible by sparse residual dependence, while a much weaker diffuse direction may not be (Example~\ref{ex:sr}).

We consider inference for the minimum dimension as such explanations approach one another. The distance to a lower-dimensional class quantifies the separation even when a particular decomposition is nonunique. Its tangent geometry gives a local first-order scale; the number and organisation of competing branches can create additional statistical difficulty. A further issue is whether observation preserves the distinguishing differences between explanations.

Sections~\ref{sec:setting}--\ref{sec:angle} collect general metric and variational tools in this setting. Sections~\ref{sec:degree}--\ref{sec:fisher} develop explicit degree-constrained prototypes and local Gaussian formulas. Section~\ref{sec:conj} formulates separated testing problems and conjectured rate benchmarks. All statements presented as results have proofs; the composite boundaries are open questions.

\paragraph{Results and organisation.}
\begin{enumerate}[leftmargin=1.5em]
\item Section~\ref{sec:setting} defines structural dimension and the margin $\Delta_q$ and shows that certified bounds on $\Delta_q$ and simultaneous bands on the distances yield honest confidence statements (Propositions~\ref{prop:cert} and~\ref{prop:rankset}).
\item Section~\ref{sec:angle} identifies the first-order behaviour of the margin with the distance of the entering direction to a tangent cone (Proposition~\ref{prop:angle}).
\item Section~\ref{sec:degree} treats competitors with bounded residual degree: the tangent cone is a union of subspaces indexed by the \emph{support slack} (Proposition~\ref{prop:slack}); distances and angles are maximum-weight $b$-matchings (Theorem~\ref{thm:bmatch}); equal-magnitude spikes have an exact angle (Theorem~\ref{thm:spike}); a shared residual support across a field creates separation absent pointwise (Proposition~\ref{prop:field}); positive semidefiniteness is bracketed by computable certificates (Proposition~\ref{prop:psd}).
\item Section~\ref{sec:fisher} gives a dimension-free local likelihood expansion for Gaussian covariance experiments (Lemma~\ref{lem:lan}), an exact formula for the Fisher information that remains after projecting out the tangent space of a rank-$r$ background (Proposition~\ref{prop:quotient}), and a detection floor in terms of the margin (Proposition~\ref{prop:floor}).
\end{enumerate}

\paragraph{Notation.}
$\Sym^p$ is the space of real symmetric $p\times p$ matrices with Frobenius inner product $\langle A,B\rangle=\tr(AB)$ and norm $\|\cdot\|_F$; $\|\cdot\|_2$ is the operator norm, $A\succeq0$ ($A\succ0$) means positive semidefinite (definite), $\lambda_{\min}$ is the smallest eigenvalue and $[x]_+=\max\{x,0\}$. For a set $C$ and a point $Y$ in a normed space, $d(Y,C)=\inf_{Z\in C}\|Y-Z\|$. Unordered pairs $\{i,j\}$, $i\ne j$, of $[p]=\{1,\dots,p\}$ are called edges; for an edge set $E$, $\deg_E(i)$ is the number of edges of $E$ containing $i$.

\section{Structural dimension and its margin}\label{sec:setting}

Let $Y\in\Sym^p$ be an observable second-order object, typically a covariance matrix, and let $\|\cdot\|$ be a norm on $\Sym^p$. A \emph{competing class} $\cS\subseteq\Sym^p$ is a closed set containing $0$, declared on scientific grounds before any rank is inferred. Fix a magnitude bound $\rho>0$ and, for $r\ge0$, let $\cL_r=\{L\in\Sym^p:L\succeq0,\ \operatorname{rank}L\le r,\ \|L\|\le\rho\}$, a compact set. The nested classes
\[
\M_{\le r}=\cL_r+\cS,\qquad \M_{\le0}\subseteq\M_{\le1}\subseteq\cdots,
\]
are closed, as sums of a compact and a closed set. Assume $Y\in\M_{\le p}$.

\begin{definition}
The \emph{structural dimension} of $Y$ is $\qstr(Y)=\min\{r:Y\in\M_{\le r}\}$. For $r\ge0$ let $D_r(Y)=d(Y,\M_{\le r})$, and for $q\ge1$ define the \emph{margin} $\Delta_q(Y)=D_{q-1}(Y)$, the amount by which the explanation with $q-1$ latent directions fails.
\end{definition}

Since the classes are closed, $\qstr(Y)=\min\{r:D_r(Y)=0\}$, and $Y$ requires a $q$-th direction exactly when $\Delta_q(Y)>0$.

\begin{remark}\label{rem:heywood}
The bound $\rho$ cannot be replaced by a closure of the unbounded sum without changing the meaning of $\qstr$. Let $\cS$ be the diagonal matrices (of either sign) and $w_n=(n,1/n,1/n)$. Then $w_nw_n^\top-\operatorname{diag}(w_nw_n^\top)$ converges to the matrix $Y$ with $Y_{12}=Y_{21}=Y_{13}=Y_{31}=1$ and all other entries $0$, so $Y$ lies in the closure of the unbounded rank-one class; yet $Y$ has no exact rank-one-plus-diagonal decomposition, which would require $w_1w_2=w_1w_3=1$ and $w_2w_3=0$. Its closure dimension $1$ would be realised only by diverging components, a Heywood-type phenomenon. If the competing class consists of positive semidefinite matrices, no bound is needed: two positive semidefinite summands of a bounded sum are bounded.
\end{remark}

\begin{lemma}\label{lem:basic}
\begin{enumerate}[label=(\roman*),leftmargin=2em]
\item $r\mapsto D_r(Y)$ is nonincreasing.
\item $|D_r(Y)-D_r(Y')|\le\|Y-Y'\|$ for all $Y,Y'$.
\item If $\M'$ is a closed set at Hausdorff distance at most $h$ from $\M_{\le r}$, then $|d(Y,\M')-D_r(Y)|\le h$.
\end{enumerate}
\end{lemma}

\begin{proof}
(i) follows from nestedness. (ii) For any $Z$, $\|Y-Z\|\le\|Y-Y'\|+\|Y'-Z\|$; take the infimum over $Z\in\M_{\le r}$ and exchange $Y,Y'$. (iii) Every $Z\in\M_{\le r}$ has a point of $\M'$ within $h$, so $d(Y,\M')\le D_r(Y)+h$, and symmetrically.
\end{proof}

Lemma~\ref{lem:basic} and Propositions~\ref{prop:cert}--\ref{prop:rankset} below use only nestedness and the triangle inequality; they hold verbatim in any normed space, provided the nested classes are closed. There closedness is an assumption: in $L^2(\mu;\Sym^p)$ bounded field classes need not be compact, so $\M_{\le r}$, $r\ge1$, need not be closed; the fields of Section~\ref{sec:degree} use only $\M_{\le0}$, which is closed. Finite dimension is used only in Proposition~\ref{prop:angle}.

\paragraph{Certificates.} The exact margin is generally the value of a nonconvex problem. A \emph{lower certificate} is a computable $\Dl_q(Z)\le\Delta_q(Z)$, typically the distance to a superset of $\M_{\le q-1}$ or a dual bound; an \emph{upper certificate} is $\Du_q(Z)=\|Z-Z_0\|$ for an explicit feasible $Z_0\in\M_{\le q-1}$.

\begin{proposition}\label{prop:cert}
Let $\widehat Y$ be an estimate of $Y$ and let $\mathcal C$ be an event on which $\|\widehat Y-Y\|\le\varepsilon$. If $\Dl_q(\widehat Y)\le\Delta_q(\widehat Y)\le\Du_q(\widehat Y)$, then on $\mathcal C$
\[
\big[\Dl_q(\widehat Y)-\varepsilon\big]_+\le\Delta_q(Y)\le\Du_q(\widehat Y)+\varepsilon .
\]
In particular, if $P(\mathcal C)\ge1-\alpha$, the interval is a confidence interval of level $1-\alpha$ for $\Delta_q(Y)$, and $\Dl_q(\widehat Y)>\varepsilon$ certifies at that level that $\qstr(Y)\ge q$.
\end{proposition}

\begin{proof}
Lemma~\ref{lem:basic}(ii) with $Y'=\widehat Y$ gives $|\Delta_q(Y)-\Delta_q(\widehat Y)|\le\varepsilon$ on $\mathcal C$; combine with the certificates and with $\Delta_q\ge0$.
\end{proof}

\begin{proposition}\label{prop:rankset}
Let $[B_r^-,B_r^+]$, $r=0,\dots,p$, be intervals with $P\big(\bigcap_{r=0}^p\{B_r^-\le D_r(Y)\le B_r^+\}\big)\ge1-\alpha$. Set $B_r^{-,\downarrow}=\max_{j\ge r}B_j^-$, $B_r^{+,\downarrow}=\min_{j\le r}B_j^+$ and $B_{-1}^{+,\downarrow}=+\infty$. For a tolerance $\tau\ge0$ let $r^*(\tau)=\min\{r:D_r(Y)\le\tau\}$ and
\[
\widehat{\mathcal Q}_\tau=\big\{r\in\{0,\dots,p\}:\ B_{r-1}^{+,\downarrow}>\tau,\ B_r^{-,\downarrow}\le\tau\big\}.
\]
Then $P(r^*(\tau)\in\widehat{\mathcal Q}_\tau)\ge1-\alpha$, and $r^*(0)=\qstr(Y)$.
\end{proposition}

\begin{proof}
On the coverage event, monotonicity (Lemma~\ref{lem:basic}(i)) gives $D_r\ge D_j\ge B_j^-$ for $j\ge r$ and $D_r\le D_j\le B_j^+$ for $j\le r$, hence $B_r^{-,\downarrow}\le D_r\le B_r^{+,\downarrow}$. For $r=r^*$, $B_{r}^{-,\downarrow}\le D_{r}\le\tau$ and, if $r\ge1$, $B_{r-1}^{+,\downarrow}\ge D_{r-1}>\tau$.
\end{proof}

\begin{remark}
Both tightened bands are nonincreasing in $r$, so $\widehat{\mathcal Q}_\tau$ is the integer interval $[r_-,r_+]$ with $r_-=\min\{r:B_r^{-,\downarrow}\le\tau\}$ and $r_+=\max\{r:B_{r-1}^{+,\downarrow}>\tau\}$, with both extrema restricted to $\{0,\ldots,p\}$ and $\min\varnothing=p+1$. The interval is empty if $r_->r_+$, which may occur outside the simultaneous coverage event. Otherwise it is a set of compatible dimensions rather than a point. The lower band gives the honest lower confidence bound $r_-$, the upper band the upper bound $r_+$. At $\tau=0$, $r_-$ is a lower confidence bound for $\qstr(Y)$, while $r_+$ is informative only when an upper certificate vanishes exactly; a positive tolerance makes both ends informative.
\end{remark}

\section{The structural angle}\label{sec:angle}

For a closed set $C$ and $x\in C$, the \emph{tangent (contingent) cone} is
\[
\Tan_C(x)=\{h:\ \exists\,t_n\downarrow0,\ h_n\to h\ \text{with}\ x+t_nh_n\in C\}
\]
\citep{AubinFrankowska1990,RockafellarWets1998}.

\begin{proposition}\label{prop:angle}
Let $C$ be closed in a finite-dimensional normed space, $x\in C$ and $h\ne0$. Then
\[
\liminf_{t\downarrow0}\frac{d(x+th,C)}{t}=d\big(h,\Tan_C(x)\big).
\]
The same holds for $d(Y_t,C)$ along any path $Y_t=x+th+o(t)$.
\end{proposition}

\begin{proof}
Write $\ell$ for the left side and $T=\Tan_C(x)$. Upper bound: $T$ is closed, so let $w\in T$ attain $d(h,T)$, with $t_n\downarrow0$, $w_n\to w$, $x+t_nw_n\in C$. Then $d(x+t_nh,C)/t_n\le\|h-w_n\|\to\|h-w\|$, so $\ell\le d(h,T)$. Lower bound: let $t_n\downarrow0$ realise $\ell$ and let $c_n\in C$ be nearest to $x+t_nh$. With $w_n=(c_n-x)/t_n$, $\|h-w_n\|\to\ell$, so $(w_n)$ is bounded and has a limit point $w$, which lies in $T$ by definition; hence $d(h,T)\le\|h-w\|=\ell$. The path statement follows from Lemma~\ref{lem:basic}(ii) applied to $C$.
\end{proof}

\begin{definition}
For $Y_0\in\M_{\le q-1}$ and a direction $h\ne0$, the \emph{structural angle} is $\kappa(Y_0,h)=d(h,\Tan_{\M_{\le q-1}}(Y_0))/\|h\|$. By Proposition~\ref{prop:angle}, $\liminf_{t\downarrow0}\Delta_q(Y_0+th)/t=\|h\|\kappa(Y_0,h)$: the margin grows at least at rate $\|h\|\kappa$ as $t\downarrow0$, and exactly at that rate along some sequence.
\end{definition}

Equality of the full limit requires regularity of $C$ at $x$. The following case covers all the examples below.

\begin{lemma}\label{lem:union}
Suppose that, in a neighbourhood $U$ of $x$, $C\cap U=\big(\bigcup_{k=1}^K(x+V_k)\big)\cap U$ for linear subspaces $V_1,\dots,V_K$. Then $\Tan_C(x)=\bigcup_kV_k$ and $d(x+th,C)=t\,d(h,\Tan_C(x))$ for all sufficiently small $t>0$.
\end{lemma}

\begin{proof}
The tangent cone is local, and that of a finite union of affine subspaces through $x$ is the union of their directions. Since $d(x+th,C)\le t\|h\|$, nearest points lie within $2t\|h\|$ of $x$, hence in $U$ for small $t$; there $d(x+th,C)=\min_kd(th,V_k)=t\min_kd(h,V_k)$.
\end{proof}

\section{Competitors with bounded residual degree}\label{sec:degree}

Throughout this section $\mathcal V=\Sym^p$ with the Frobenius norm; the working-norm normalisation $p^{-1}\|\cdot\|_F$ multiplies all distances by $p^{-1}$ and leaves angles unchanged. For an edge set $E$ let
\[
\Omega_E=\{S\in\Sym^p:\ S_{ij}=0\ \text{for}\ i\ne j,\ \{i,j\}\notin E\},
\]
a linear subspace with free diagonal. For a degree budget $\dbar\in\mathbb Z_{\ge0}$ the competing class of residuals with at most $\dbar$ off-diagonal entries per row is
\[
\cS_{\dbar}=\bigcup_{E\in\cB_{\dbar}}\Omega_E,\qquad \cB_{\dbar}=\{E:\ \deg_E(i)\le\dbar\ \text{for all}\ i\},
\]
a closed cone (a finite union of subspaces). Magnitude bounds on the residual entries are omitted: all local statements below concern points at which such bounds are slack. With this competitor $\M_{\le0}=\cS_{\dbar}$.

\begin{lemma}\label{lem:proj}
For $Y\in\Sym^p$, $d_F(Y,\Omega_E)^2=2\sum_{\{i,j\}\notin E}Y_{ij}^2$, attained by the matrix that copies $Y$ on the diagonal and on $E$ and vanishes elsewhere.
\end{lemma}

\begin{proof}
$\Omega_E$ is a coordinate subspace of $\Sym^p$; orthogonal projection keeps the free coordinates. Each off-diagonal pair contributes twice to the Frobenius norm.
\end{proof}

\subsection{\texorpdfstring{Distances as $b$-matchings}{Distances as b-matchings}}

Given nonnegative edge weights $w$ and vertex capacities $b_i\in\mathbb Z_{\ge0}$, a \emph{maximum-weight $b$-matching} is an edge set $E$ maximising $\sum_{e\in E}w_e$ subject to $\deg_E(i)\le b_i$. It is solvable in polynomial time \citep[Chapters~31--33]{Schrijver2003}.

\begin{theorem}\label{thm:bmatch}
For $Y\in\Sym^p$,
\[
d_F(Y,\cS_{\dbar})^2=W(Y)-\max_{E\in\cB_{\dbar}}\sum_{\{i,j\}\in E}w_{ij}(Y),\qquad w_{ij}(Y)=2Y_{ij}^2,\quad W(Y)=\sum_{i<j}w_{ij}(Y),
\]
so the rank-zero margin $\Delta_1(Y)$ is computable exactly in polynomial time. In the working norm, $\Delta_1(Y)^2=W-\max_{E\in\cB_{\dbar}}\sum_Ew_{ij}$ with $w_{ij}=2p^{-2}Y_{ij}^2$.
\end{theorem}

\begin{proof}
$d_F(Y,\cS_{\dbar})^2=\min_{E\in\cB_{\dbar}}d_F(Y,\Omega_E)^2$, and by Lemma~\ref{lem:proj} this equals $\min_E\big(W-\sum_Ew_{ij}\big)$.
\end{proof}

\begin{example}[strong but redundant, weak but necessary]\label{ex:sr}
Let $\dbar\ge1$. For $v=(e_1+e_2)/\sqrt2$ and any $a>0$, $Y=a\,vv^\top$ has a single nonzero off-diagonal pair, which a one-edge $E$ covers: $\Delta_1(Y)=0$ however large the eigenvalue $a$. For $v=p^{-1/2}\mathbf 1$ and $a>0$, all $p(p-1)/2$ pairs carry weight $2a^2/p^2$ and at most $\lfloor p\dbar/2\rfloor$ can be covered, so $d_F(Y,\cS_{\dbar})^2\ge a^2(1-\tfrac1p-\tfrac{\dbar}{p})>0$ for $\dbar<p-1$, however small $a$ is.
\end{example}

\subsection{Tangent cones and support slack}

Let $S_0\in\cS_{\dbar}$ have active off-diagonal support $E_0=\{\{i,j\}:(S_0)_{ij}\ne0\}$ and define the \emph{support slack} $c_i=\dbar-\deg_{E_0}(i)\ge0$ and
\[
\cE(E_0,c)=\{E'\ \text{disjoint from}\ E_0:\ \deg_{E'}(i)\le c_i\ \text{for all}\ i\}.
\]

\begin{proposition}\label{prop:slack}
$\Tan_{\cS_{\dbar}}(S_0)=\bigcup_{E'\in\cE(E_0,c)}\Omega_{E_0\cup E'}$, and for every $h$, $d_F(S_0+th,\cS_{\dbar})=t\,d_F(h,\Tan_{\cS_{\dbar}}(S_0))$ for small $t>0$.
\end{proposition}

\begin{proof}
Let $\delta=\min_{E_0}|(S_0)_{ij}|>0$ (or $\delta=\infty$ if $E_0=\emptyset$). If $\max_{ij}|S_{ij}-(S_0)_{ij}|<\delta$, the off-diagonal support of $S$ contains $E_0$, so $S\in\cS_{\dbar}$ if and only if its support is $E_0\cup E'$ with $E'\in\cE(E_0,c)$, i.e.\ $S\in\Omega_{E_0\cup E'}$. Conversely each $\Omega_{E_0\cup E'}\subseteq\cS_{\dbar}$ and contains $S_0$. Apply Lemma~\ref{lem:union}.
\end{proof}

Saturated supports ($c\equiv0$) are regular points with a single tangent subspace; an empty support ($c\equiv\dbar$) is maximally singular, with tangent cone the whole of $\cS_{\dbar}$.

\begin{corollary}\label{cor:energy}
Let $h=vv^\top$ with $\|v\|_2=1$, so $\|h\|_F=1$. Then
\[
\kappa(S_0,h)^2=1-\sum_iv_i^4-2\!\!\sum_{\{i,j\}\in E_0}\!\!v_i^2v_j^2-\max_{E'\in\cE(E_0,c)}2\!\!\sum_{\{i,j\}\in E'}\!\!v_i^2v_j^2,
\]
the last term being a maximum-weight $b$-matching with capacities $b_i=c_i$ on the complement of $E_0$.
\end{corollary}

\begin{proof}
Combine Proposition~\ref{prop:slack} with Lemma~\ref{lem:proj} and $2\sum_{i<j}v_i^2v_j^2=(\sum_iv_i^2)^2-\sum_iv_i^4=1-\sum_iv_i^4$.
\end{proof}

\subsection{The exact angle of equal-magnitude spikes}

\begin{lemma}\label{lem:factor}
For integers $k\ge1$ and $0\le d\le k-1$, the complete graph $K_k$ has a subgraph of maximum degree at most $d$ with exactly $\lfloor kd/2\rfloor$ edges, and no such subgraph has more.
\end{lemma}

\begin{proof}
The upper bound is the handshake identity $2|E|=\sum_i\deg_E(i)\le kd$. For the construction label the vertices by $\mathbb Z_k$ and let $G_s$ be the circulant graph with edges $\{i,i\pm1\},\dots,\{i,i\pm s\}$. If $d=k-1$ take $K_k$. If $d<k-1$ is even, $G_{d/2}$ is $d$-regular since $d/2<k/2$. If $d$ is odd and $k$ even, add to $G_{(d-1)/2}$ the perfect matching $\{i,i+k/2\}$, which is disjoint from it because $(d-1)/2<k/2$; the result is $d$-regular. If $d$ is odd and $k$ odd, then $d\le k-2$; the edges $\{i,i+(k-1)/2\}$ are not in $G_{(d-1)/2}$ (their difference $\pm(k-1)/2$ exceeds $(d-1)/2$ in absolute value modulo $k$) and form a Hamiltonian cycle, because $\gcd((k-1)/2,k)=1$; it contains a matching with $(k-1)/2$ edges, and adding it to $G_{(d-1)/2}$ gives $k(d-1)/2+(k-1)/2=\lfloor kd/2\rfloor$ edges with maximum degree $d$.
\end{proof}

\begin{theorem}\label{thm:spike}
Let $v$ have support $K\subseteq[p]$ with $|K|=k$ and $|v_i|=k^{-1/2}$ on $K$, let $0\le\dbar\le k-1$ and $g=k-1-\dbar$. At a residual $S_0$ with empty off-diagonal support (for instance $S_0=0$ or diagonal), the structural angle of $h=vv^\top$ with respect to $\M_{\le0}=\cS_{\dbar}$ is
\[
\kappa_{k,\dbar}^2=\frac{2\lceil kg/2\rceil}{k^2}=\frac gk+O(k^{-2}).
\]
\end{theorem}

\begin{proof}
By Corollary~\ref{cor:energy} with $E_0=\emptyset$ and $c\equiv\dbar$, the weights $2v_i^2v_j^2$ equal $2/k^2$ on pairs inside $K$ and vanish elsewhere, and $1-\sum_iv_i^4=1-1/k=(2/k^2)\binom k2$. By Lemma~\ref{lem:factor} at most $\lfloor k\dbar/2\rfloor$ pairs inside $K$ can be covered, and this is attained. Hence $\kappa^2=(2/k^2)\big(\binom k2-\lfloor k\dbar/2\rfloor\big)$. Since $k(k-1)$ is even, $kg=k(k-1)-k\dbar$ has the parity of $k\dbar$, and $\binom k2-\lfloor k\dbar/2\rfloor=\lceil kg/2\rceil$.
\end{proof}

\begin{proposition}[equal angle, different branching]\label{prop:clique}
Let $k\ge2$ be even, $\dbar=k-2$ and $v$ as in Theorem~\ref{thm:spike}.
\begin{enumerate}[label=(\roman*),leftmargin=2em]
\item At an empty support, $\kappa^2=1/k$, and the minimising tangent subspaces are exactly the $\Omega_{E}$ whose restriction $E\cap\binom K2$ is the complement of a perfect matching of $K$: there are $(k-1)!!=(k-1)(k-3)\cdots1$ such restrictions.
\item At a residual $S_0$ whose support on $K$ is the complete graph minus a fixed perfect matching $\mathcal P$ (every vertex of $K$ saturated), $\kappa^2=1/k$ and the tangent cone restricted to $K$ is the single subspace $\Omega_{E_0}$.
\end{enumerate}
\end{proposition}

\begin{proof}
(i) Theorem~\ref{thm:spike} with $g=1$ gives $\kappa^2=2(k/2)/k^2=1/k$. An optimal $E$ inside $K$ has $k(k-2)/2$ edges and maximum degree $k-2$, hence is $(k-2)$-regular on $K$, and its complement in $K_k$ is $1$-regular, i.e.\ a perfect matching; there are $(k-1)!!$ perfect matchings of $K_k$. (ii) The slack vanishes on $K$, so by Proposition~\ref{prop:slack} no edge inside $K$ can be added, and the uncovered pairs are those of $\mathcal P$: $\kappa^2=2\cdot\frac k2\cdot\frac1{k^2}=1/k$.
\end{proof}

The two backgrounds share degree budget and angle; they differ in the number of nuisance branches through which the spike can be explained. This is the reason for Conjecture~\ref{conj:static}.

\subsection{Fields with a shared residual support}

In this subsection the ambient space is $L^2(\mu;\Sym^p)$, possibly infinite-dimensional; only $\M_{\le0}$ and Lemma~\ref{lem:basic} are used. Let $Y(u)\in\Sym^p$, $u\in\mathcal U$, be measurable and square-integrable with respect to a probability measure $\mu$, with norm $\|Y\|^2=\int\|Y(u)\|_F^2\,d\mu(u)$. Compare the \emph{pointwise} competitor, in which $S(u)\in\cS_{\dbar}$ for $\mu$-almost every $u$ with a support that may depend on $u$, with the \emph{shared-support} competitor $\cS^{\mathrm{sh}}_{\dbar}=\bigcup_{E\in\cB_{\dbar}}\{S:\ S(u)\in\Omega_E\ \text{for }\mu\text{-almost every }u\}$.

\begin{proposition}\label{prop:field}
With $\bar w_{ij}=2\int Y_{ij}(u)^2\,d\mu(u)$ and $\bar W=\sum_{i<j}\bar w_{ij}$,
\[
d(Y,\cS^{\mathrm{sh}}_{\dbar})^2=\bar W-\max_{E\in\cB_{\dbar}}\sum_{\{i,j\}\in E}\bar w_{ij},
\qquad
d(Y,\text{pointwise})^2=\int d_F\big(Y(u),\cS_{\dbar}\big)^2d\mu(u)\le d(Y,\cS^{\mathrm{sh}}_{\dbar})^2,
\]
and the inequality can be strict with the left side equal to $0$.
\end{proposition}

\begin{proof}
For fixed $E$ the projection acts pointwise (Lemma~\ref{lem:proj}), giving $\int2\sum_{\{i,j\}\notin E}Y_{ij}^2\,d\mu=\bar W-\sum_E\bar w_{ij}$; minimise over $E$. The pointwise problem decouples across $u$ (a measurable selection of optimal supports exists because $\cB_{\dbar}$ is finite), and the shared class is contained in the pointwise one. Example: $p=3$, $\dbar=1$, $\mu$ uniform on $\{u_1,u_2\}$, $Y(u_1)$ supported on the pair $\{1,2\}$ and $Y(u_2)$ on $\{1,3\}$ with unit entries. Pointwise each matrix is covered by one edge, so the pointwise distance is $0$; a shared $E$ with $\deg_E(1)\le1$ cannot contain both pairs, so $d(Y,\cS^{\mathrm{sh}}_1)^2=\min(\bar w_{12},\bar w_{13})=1>0$.
\end{proof}

Positive structural dimension can therefore be certified from the whole field even when every pointwise problem is explained by the competitor.

\subsection{Certifying positive semidefiniteness}

Covariance residuals are positive semidefinite. Let $\cS^+_{\dbar}=\cS_{\dbar}\cap\{S\succeq0\}$, and for a fixed $E$ write $D_{\mathrm{rel},E}=d_F(Y,\Omega_E)$, $D_{\mathrm{PSD},E}=d_F(Y,\Omega_E\cap\{S\succeq0\})$, $S^{\mathrm{rel}}_E$ for the projection of Lemma~\ref{lem:proj} and $\zeta_E=[-\lambda_{\min}(S^{\mathrm{rel}}_E)]_+$.

\begin{proposition}\label{prop:psd}
$D_{\mathrm{rel},E}\le D_{\mathrm{PSD},E}\le\big(D_{\mathrm{rel},E}^2+p\,\zeta_E^2\big)^{1/2}$. Consequently
\[
d_F(Y,\cS_{\dbar})\ \le\ d_F(Y,\cS^+_{\dbar})\ \le\ \min_{E\in\cB_{\dbar}}\big(D_{\mathrm{rel},E}^2+p\,\zeta_E^2\big)^{1/2},
\]
where the left side is the $b$-matching value of Theorem~\ref{thm:bmatch} and the right side can be evaluated at any computed support, for instance the $b$-matching optimiser. In the working norm the correction term becomes $p^{-1}\zeta_E^2$.
\end{proposition}

\begin{proof}
The lower bounds hold because the PSD classes are subsets. $S^{\mathrm{rel}}_E+\zeta_EI$ is positive semidefinite, belongs to $\Omega_E$ and differs from $S^{\mathrm{rel}}_E$ only on the diagonal, where $Y-S^{\mathrm{rel}}_E$ vanishes; hence $\|Y-S^{\mathrm{rel}}_E-\zeta_EI\|_F^2=D_{\mathrm{rel},E}^2+p\,\zeta_E^2$. Multiplying distances by $p^{-1}$ turns $p\zeta_E^2$ into $p^{-1}\zeta_E^2$.
\end{proof}

The same argument shows $D_{\mathrm{PSD},E}^2\le D_{\mathrm{rel},E}^2+\|\Lambda\|_F^2$ for any diagonal $\Lambda$ with $S^{\mathrm{rel}}_E+\Lambda\succeq0$. The scalar shift is the simplest choice and can be conservative when the eigenvalue deficit is localised on few coordinates; minimising $\|\Lambda\|_F$ under this constraint is a semidefinite programme.

\section{Local Fisher structure}\label{sec:fisher}

Let $z_1,\dots,z_T$ be independent $N(0,\Sigma)$ vectors in $\mathbb R^p$ and $P_\Sigma^{\otimes T}$ their joint law. The Fisher inner product at $\Sigma_0\succ0$ is
\[
\langle A,B\rangle_{I,0}=\tfrac12\tr\big(\Sigma_0^{-1}A\Sigma_0^{-1}B\big),\qquad \|A\|_{I,0}^2=\langle A,A\rangle_{I,0}.
\]

\begin{lemma}[dimension-free likelihood expansion]\label{lem:lan}
Let $\Sigma_0=I$ and let $A=A_{p,T}\in\Sym^p$ satisfy $T\|A\|_F^2=O(1)$ (hence $\|A\|_2\le\|A\|_F\to0$), with $p=p_T$ arbitrary. Under $P_I^{\otimes T}$,
\[
\log\frac{dP_{I+A}^{\otimes T}}{dP_I^{\otimes T}}=\frac12\sum_{t=1}^T\big(z_t^\top Az_t-\tr A\big)-\frac T4\tr(A^2)+o_P(1).
\]
If moreover $\frac T2\tr(A^2)\to\sigma^2$, the right side converges in law to $N(-\sigma^2/2,\sigma^2)$, and $P_{I+A}^{\otimes T}$ and $P_I^{\otimes T}$ are mutually contiguous.
\end{lemma}

\begin{proof}
Assume $\|A\|_2\le1/2$ and let $R=(I+A)^{-1}$, so $\|R\|_2\le2$ and $R-I=-A+A^2R$. The log-likelihood ratio is
\[
\ell_T=-\frac T2\log\det(I+A)+\frac12\sum_t z_t^\top Az_t-\frac12\sum_tz_t^\top A^2Rz_t .
\]
For $|\lambda|\le1/2$, $|\log(1+\lambda)-\lambda+\lambda^2/2|\le|\lambda|^3$, so $\log\det(I+A)=\tr A-\frac12\tr A^2+r_1$ with $|r_1|\le\|A\|_2\|A\|_F^2$, and $T|r_1|\to0$. The matrix $A^2R$ is symmetric with $\|A^2R\|_F\le2\|A\|_2\|A\|_F$; hence $\sum_tz_t^\top A^2Rz_t$ has mean $T\tr(A^2R)=T\tr A^2-T\tr(A^3R)$, with $T|\tr(A^3R)|\le2T\|A\|_2\|A\|_F^2\to0$, and variance $2T\|A^2R\|_F^2\le8\|A\|_2^2\,T\|A\|_F^2\to0$. Collecting terms gives the expansion. For the limit law, $X_t=\frac12(z_t^\top Az_t-\tr A)$ are i.i.d., centred, with variance $\frac12\tr A^2$; as second-order Gaussian chaoses they satisfy $E|X_t|^3\le C(\tr A^2)^{3/2}$ by hypercontractivity \citep{Janson1997}, so for $\sigma^2>0$ the Lyapunov ratio is $O(T^{-1/2})$ and $\sum_tX_t\Rightarrow N(0,\sigma^2)$; for $\sigma^2=0$ the variance $\frac T2\tr A^2$ of $\sum_tX_t$ tends to $0$. Contiguity follows from Le Cam's first lemma \citep{LeCamYang2000}.
\end{proof}

For general $\Sigma_0$, whitening with $\Sigma_0^{-1/2}$ maps a perturbation $G$ to $A=\Sigma_0^{-1/2}G\Sigma_0^{-1/2}$, and $\frac T2\tr(A^2)=T\|G\|_{I,0}^2$. In particular a perturbation is locally detectable at scale $(\sqrt T\|G\|_{I,0})^{-1}$.

\begin{corollary}\label{cor:kl}
If $M=\Sigma^{-1/2}G\Sigma^{-1/2}$ satisfies $\|M\|_2\to0$, then $\mathrm{KL}\big(P_\Sigma^{\otimes T}\,\|\,P_{\Sigma+G}^{\otimes T}\big)=\frac T4\|M\|_F^2\{1+O(\|M\|_2)\}=\frac T2\|G\|_{I,\Sigma}^2\{1+O(\|M\|_2)\}$.
\end{corollary}

\begin{proof}
With eigenvalues $\lambda_j$ of $M$, $\mathrm{KL}=\frac T2\sum_j\phi(\lambda_j)$ where $\phi(\lambda)=\frac1{1+\lambda}-1+\log(1+\lambda)=\frac{\lambda^2}2+O(|\lambda|^3)$.
\end{proof}

\subsection{\texorpdfstring{The Fisher quotient at a rank-$r$ background}{The Fisher quotient at a rank-r background}}

When the null model contains a rank-$r$ component, a perturbation is informative only through the part not absorbed by moving that component. Let $\Sigma_0=L_0+S_0$ with $L_0=U\Lambda U^\top$, $U\in\mathbb R^{p\times r}$ with orthonormal columns, $\Lambda\succ0$ diagonal and $S_0\succ0$, and let $U_\perp$ have orthonormal columns spanning the orthogonal complement of the range of $U$. Assume $\|L_0\|<\rho$, so that the magnitude bound of Section~\ref{sec:setting} is slack; then the tangent cone of $\cL_r$ at $L_0$ is the tangent space of the rank-$r$ positive semidefinite matrices,
\[
\Tan_U=\{UB^\top+BU^\top:\ B\in\mathbb R^{p\times r}\}=\{Q\in\Sym^p:\ U_\perp^\top QU_\perp=0\}
\]
\citep{Absil2008}.

\begin{proposition}\label{prop:quotient}
For $G\in\Sym^p$, with $G_{22}=U_\perp^\top GU_\perp$ and $S_\perp=U_\perp^\top S_0U_\perp\succ0$,
\[
\min_{Q\in\Tan_U}\|G-Q\|_{I,0}^2=\tfrac12\big\|S_\perp^{-1/2}G_{22}S_\perp^{-1/2}\big\|_F^2 .
\]
In particular the Fisher information left after projecting out $\Tan_U$ does not depend on the strength $\Lambda$ of the background.
\end{proposition}

\begin{proof}
Write $M=G-Q$. As $Q$ ranges over $\Tan_U$, $M$ ranges over $\{M\in\Sym^p:U_\perp^\top MU_\perp=G_{22}\}$, so we minimise the convex function $\frac12\tr(\Sigma_0^{-1}M\Sigma_0^{-1}M)$ under a linear constraint. Stationarity gives $\Sigma_0^{-1}M\Sigma_0^{-1}=U_\perp\Gamma U_\perp^\top$ for a symmetric multiplier $\Gamma$, i.e.\ $M=\Sigma_0U_\perp\Gamma U_\perp^\top\Sigma_0$. Since $U_\perp^\top L_0=0$, $U_\perp^\top\Sigma_0U_\perp=S_\perp$, and the constraint becomes $S_\perp\Gamma S_\perp=G_{22}$, so $\Gamma=S_\perp^{-1}G_{22}S_\perp^{-1}$. The value is $\frac12\tr(U_\perp\Gamma U_\perp^\top\Sigma_0U_\perp\Gamma U_\perp^\top\Sigma_0)=\frac12\tr(\Gamma S_\perp\Gamma S_\perp)=\frac12\tr(S_\perp^{-1}G_{22}S_\perp^{-1}G_{22})$.
\end{proof}

\begin{lemma}\label{lem:coord}
Let $G$ be supported on a coordinate set $\mathcal I$, i.e.\ $G=E_{\mathcal I}G_{\mathcal I}E_{\mathcal I}^\top$ with $E_{\mathcal I}$ the coordinate embedding, and let $\vartheta=\|U^\top E_{\mathcal I}\|_2^2\le1$. Then
\[
(1-\vartheta)\|G\|_F\le\|U_\perp^\top GU_\perp\|_F\le\|G\|_F,
\]
and therefore $\frac{(1-\vartheta)^2}{2\|S_\perp\|_2^2}\|G\|_F^2\le\min_{Q\in\Tan_U}\|G-Q\|_{I,0}^2\le\frac{1}{2\lambda_{\min}(S_\perp)^2}\|G\|_F^2$.
\end{lemma}

\begin{proof}
Let $N=E_{\mathcal I}^\top U_\perp U_\perp^\top E_{\mathcal I}=I-E_{\mathcal I}^\top UU^\top E_{\mathcal I}$, whose eigenvalues lie in $[1-\vartheta,1]$. Then $\|U_\perp^\top GU_\perp\|_F^2=\tr(G_{\mathcal I}NG_{\mathcal I}N)=\|N^{1/2}G_{\mathcal I}N^{1/2}\|_F^2\in[(1-\vartheta)^2,1]\cdot\|G_{\mathcal I}\|_F^2$. The last display follows from Proposition~\ref{prop:quotient} and $\|S_\perp^{-1/2}CS_\perp^{-1/2}\|_F\in[\|C\|_F/\|S_\perp\|_2,\ \|C\|_F/\lambda_{\min}(S_\perp)]$.
\end{proof}

For coordinate-supported perturbations, Lemma~\ref{lem:coord} controls the Frobenius norm retained after compression onto the orthogonal complement of the latent space. The corresponding Fisher information is weighted by $S_\perp^{-1}$, as specified in Proposition~\ref{prop:quotient}. For perturbations through the fixed decomposition $\Sigma_0=L_0+S_0$, allowing the residual to move along a branch $\Omega_E$ gives the Fisher distance to $\Tan_U+\Omega_E$, minimised over the branches through $S_0$. The sandwich $S_\perp^{-1/2}U_\perp^\top(\cdot)U_\perp S_\perp^{-1/2}$ couples the residual coordinate subspaces. Identifying this branchwise profile with the tangent geometry of the full structural null also requires control of other decompositions and of the projection. These joint questions remain open.

\subsection{A detection floor from the margin}

The margin bounds from below what any test can detect. Let $\cN\subseteq\{\Sigma\succ0\}$ be a null family and, for a test $\psi$, write $\mathcal R(\psi;\Sigma_0,\Sigma_1)=E_{\Sigma_0}\psi+E_{\Sigma_1}(1-\psi)$ for the sum of its two error probabilities under $P_{\Sigma_0}^{\otimes T}$ and $P_{\Sigma_1}^{\otimes T}$.

\begin{proposition}[margin floor]\label{prop:floor}
Let $\Sigma_1\succ0$ and let $\Sigma_0\in\cN$ satisfy $\lambda_{\min}(\Sigma_0)\ge c>0$ and $\|\Sigma_1-\Sigma_0\|_F\le c/2$. Then for every test $\psi$,
\[
\mathcal R(\psi;\Sigma_0,\Sigma_1)\ \ge\ 1-\frac{\sqrt T\,\|\Sigma_1-\Sigma_0\|_F}{\sqrt2\,c} .
\] In particular, let $\cN=\cS_{\dbar}\cap\{\Sigma\succ0\}$ and measure $\Delta_1$ in the unnormalised Frobenius norm. If a nearest point $\Sigma_0$ of $\cS_{\dbar}$ to $\Sigma_1$ satisfies $\lambda_{\min}(\Sigma_0)\ge c$ and $\Delta_1(\Sigma_1)\le c/2$, the bound reads $1-\sqrt{T/2}\,\Delta_1(\Sigma_1)/c$: no test separates $\Sigma_1$ from the null once $\sqrt T\,\Delta_1(\Sigma_1)\to0$, whatever $p$.
\end{proposition}

\begin{proof}
Let $G=\Sigma_1-\Sigma_0$ and $M=\Sigma_0^{-1/2}G\Sigma_0^{-1/2}$, so $\|M\|_2\le\|M\|_F\le\|G\|_F/c\le\frac12$. For $|\lambda|\le\frac12$ the function $\phi$ in the proof of Corollary~\ref{cor:kl} satisfies $0\le\phi(\lambda)\le2\lambda^2$, hence $\mathrm{KL}(P_{\Sigma_0}^{\otimes T}\|P_{\Sigma_1}^{\otimes T})\le T\|M\|_F^2\le T\|G\|_F^2/c^2$. By Pinsker's inequality the total variation distance is at most $\sqrt{T/2}\,\|G\|_F/c$, and $\mathcal R(\psi;\Sigma_0,\Sigma_1)\ge1-\mathrm{TV}$ for every test.
\end{proof}

\begin{corollary}\label{cor:floorspike}
In the setting of Proposition~\ref{prop:clique}, let $S_0\in\cS_{k-2}$ be either residual of parts (i)--(ii), with $\lambda_{\min}(S_0)\ge c>0$, and let $\Sigma_1=S_0+\theta vv^\top$ with $\theta>0$. If $\theta/\sqrt k\le c/4$, there is $\Sigma_0\in\cS_{k-2}$, $\Sigma_0\succ0$, such that every test satisfies
\[
\mathcal R(\psi;\Sigma_0,\Sigma_1)\ \ge\ 1-\frac{\sqrt2}{c}\,\theta\sqrt{T/k}.
\]
In particular no test separates $\Sigma_1$ from $\cS_{k-2}$ when $\theta\sqrt{T/k}\to0$, uniformly in $p$ and $k$ for fixed $c$.
\end{corollary}

\begin{proof}
Let $E$ be the support of a minimising subspace in Proposition~\ref{prop:clique}: in case (i) any $E$ whose restriction to $K$ is the complement of a perfect matching $\mathcal P$ of $K$, in case (ii) $E=E_0$ with the matching $\mathcal P$ of that statement. Let $\Sigma_0$ be the orthogonal projection of $\Sigma_1$ onto $\Omega_E$; since $S_0\in\Omega_E$ in both cases, $\Sigma_0\in\Omega_E\subseteq\cS_{k-2}$ and $G=\Sigma_1-\Sigma_0$ is supported on the pairs of $\mathcal P$, with entries $\pm\theta/k$. Hence $\|G\|_F=\theta/\sqrt k$ and $\|G\|_2=\theta/k$, for every $\theta>0$. By Weyl's inequality $\lambda_{\min}(\Sigma_0)\ge\lambda_{\min}(\Sigma_1)-\theta/k\ge c-c/4\ge c/2$, using $\lambda_{\min}(\Sigma_1)\ge\lambda_{\min}(S_0)$. Apply Proposition~\ref{prop:floor} with $c/2$ in place of $c$, noting $\|G\|_F\le c/4$.
\end{proof}

For small $\theta$, $\|G\|_F=\theta/\sqrt k=\theta\kappa$ is exactly the margin $\Delta_1(\Sigma_1)$ (Propositions~\ref{prop:slack} and~\ref{prop:clique}); the construction above avoids any condition on the size of the active entries of $S_0$. The floor depends on the background only through the structural angle and is therefore the same at the empty and at the saturated support. It is a single-branch bound: it ignores that, at an empty support, the null can be approached through any of $(k-1)!!$ branches.

\begin{remark}[two scales]\label{rem:scales}
Proposition~\ref{prop:cert} resolves margins above the estimation error $\varepsilon\asymp\|\widehat\Sigma-\Sigma\|_F$, which for the sample covariance with $\tr\Sigma\asymp p$ is of order $p/\sqrt T$ in the unnormalised Frobenius norm, since $E\|\widehat\Sigma-\Sigma\|_F^2=\{(\tr\Sigma)^2+\|\Sigma\|_F^2\}/T$. Proposition~\ref{prop:floor} places the information floor at $\Delta_1\asymp1/\sqrt T$. Lipschitz certificates are therefore valid but may be conservative by a factor up to $p$; certificates that exploit the structure of the competitor rather than the global norm of the estimation error are needed to close this gap, and are not treated here.
\end{remark}

\section{Conjectures and open problems}\label{sec:conj}

The degree-constrained distance and tangent formulas above account exactly for unions of support branches. The Gaussian bounds, however, do not determine the sharp testing cost of those unions. We specify the structural testing target before stating the rate predictions.

\subsection{Separated hypotheses and local benchmarks}
Let $0<c<C<\infty$ be fixed, and restrict the covariance classes to $cI\preceq\Sigma\preceq CI$. For a declared neighbourhood $\mathcal U_T$ and margin $\delta>0$, define
\[
\mathcal H_0=\M_{\le r}\cap\mathcal U_T,\qquad
\mathcal H_1(\delta)=\{\Sigma\in\M_{\le r+1}\cap\mathcal U_T:
d_F(\Sigma,\M_{\le r})\ge\delta\}.
\]
The classes here use positive-semidefinite residuals, the degree budget and any magnitude caps specified for the experiment. The margin is measured in the unnormalised Frobenius norm. For a test $\psi$ with values in $[0,1]$, its minimax sum of errors is
\[
R_T(\delta)=\inf_\psi\left\{
\sup_{\Sigma\in\mathcal H_0}E_\Sigma\psi+
\sup_{\Sigma\in\mathcal H_1(\delta)}E_\Sigma(1-\psi)\right\}.
\]
This excludes a covariance that has both a spike representation and a lower-rank representation from the alternatives. A nonzero spike coefficient alone would not ensure structural separation.

The rank-zero amplitude benchmarks use $\Sigma=S+\theta vv^\top$, with $S\in\cS_{\dbar}$ positive definite and $v$ an equal-magnitude unit vector on $k$ coordinates. For an amplitude parameter $a>0$, restrict $a\le\theta\le2a$ and require
\[
c_*\theta\kappa_{k,\dbar}\le d_F(\Sigma,\M_{\le0})
\le C_*\theta\kappa_{k,\dbar},\qquad 0<c_*<1<C_*<\infty.
\]
Thus amplitude and structural margin are linked on the declared alternatives. Both hypotheses are localised around the specified residual background, with uniform spectral bounds. In a saturated-support experiment the neighbourhood preserves the active edges; at an empty background it permits all degree-admissible support branches. Neighbourhood widths must contain the signal-scale perturbations. The precise ranges of localisation and growth for which sharp matching is possible form part of the rate problem.

For these amplitude-indexed families, a boundary $a_T$ means that the minimax error tends to zero when $a/a_T\to\infty$, and stays away from zero when $a/a_T\to0$, within the declared local regime. Unknown-direction experiments take the union over the permitted sparse directions and their corresponding separated alternatives.

\subsection{Static and spectral predictions}
\begin{conjecture}[static branching benchmark]\label{conj:static}
Consider the separated rank-zero local experiments above with a known balanced signed direction on an even number $k$ of coordinates and $\dbar=k-2$. Under suitable localisation and growth conditions and up to logarithmic factors, the boundary is $a_T\asymp k/\sqrt T$ at an empty residual support and $a_T\asymp\sqrt{k/T}$ at a saturated clique-minus-matching support: branch uncertainty costs a factor of order $\sqrt k$ despite the equal structural angles of Proposition~\ref{prop:clique}. For an unknown direction and general $0\le\dbar\le k-2$, writing $g=k-1-\dbar$, the boundary is of order
\[
\max\left\{\sqrt{\frac{k\log(ep/k)}T},\sqrt{\frac{k(\dbar+1)}{gT}}\right\}.
\]
\end{conjecture}

Corollary~\ref{cor:floorspike} establishes a $\sqrt{k/T}$ floor in both endpoint configurations. It does not prove the larger branching rate, its achievability, or uniformity over the local nuisance families. Establishing or revising the proposed orders requires composite lower bounds and matching tests; the count of branches alone does not determine their statistical cost.

For the dynamic benchmark, the data are $T$ consecutive observations of a real stationary Gaussian process. The spectral classes are Hermitian, conjugate-symmetric and uniformly bounded between $cI$ and $CI$. On $[0,\pi]$ the partition into $m$ equal bands is known; the residual support may be independently reconfigured on those bands and is reflected to negative frequencies. For $\mathcal N_m$ the corresponding rank-zero spectral class, measure separation by
\[
d_{F,2}(F,\mathcal N_m)=\inf_{H\in\mathcal N_m}
\left\{\frac1\pi\int_0^\pi\|F(\omega)-H(\omega)\|_F^2d\omega\right\}^{1/2}.
\]
The alternatives $F(\omega)=S(\omega)+\theta vv^\top$ have $a\le\theta\le2a$ and margin between $c_*\theta/\sqrt k$ and $C_*\theta/\sqrt k$. The minimax risk is defined as above using the exact finite-sample Gaussian law with block Toeplitz covariance. Unknown band locations and globally smooth support transitions define additional experiments.

\begin{conjecture}[switching benchmark]\label{conj:dyn}
For the separated stationary-Gaussian benchmark with known balanced direction on an even number $k$ of coordinates, $\dbar=k-2$, and $m=J+1$ equal bands with independently reconfigured perfect matchings, the boundary is, up to logarithmic factors,
\[
a_T\asymp\frac{k m^{1/4}}{\sqrt T},\qquad
\frac{k^2\sqrt m}{T}\to0,\quad \frac{k m^{5/2}\log(eT)}{T}\to0,
\]
in the exact Toeplitz experiment, under suitable localisation conditions.
\end{conjecture}

The fourth-root penalty is a prediction to be tested, not an implication of the deterministic switching geometry. A Whittle approximation, if used, must be justified uniformly over the declared bands and nuisance alternatives; the independent-observation expansion of Lemma~\ref{lem:lan} does not itself provide that justification.

Further open problems are joint profiling of the low-rank and residual tangents, separation bounds for growing rank, information loss under discretisation and spectral smoothing, and computational certificates narrow enough at the statistical resolution. A characterization that changes a conjectured exponent or identifies an additional necessary condition resolves part of these questions.

\end{document}